\documentstyle[12pt,twoside]{article}
\input xy
\xyoption{all}
 \newtheorem{lem}{Lemma}[section]
\newtheorem{remark}[lem]{Remark}

\newtheorem{exemple}[lem]{Example}

\newtheorem{exemples}[lem]{Examples}

\def\L{{\cal L}}

\def\A{{\cal A}}

\def\X{{\cal X}}
\def\O{{\cal O}}

\def\F{{\cal F}}

 1
\font\tenbbb=msbm10 scaled\magstep 1
\font\sevenbbb=msbm7
\font\fivebbb=msbm5
\newfam\bbbfam
\textfont\bbbfam=\tenbbb
\scriptfont\bbbfam=\sevenbbb
\scriptscriptfont\bbbfam=\fivebbb

\font\pt=cmr8
\begin{document}
\vbox to0.5truecm {\vfill}

\vglue3truecm

\begin{center}
{\Large\bf
A NOTE ON A CONJECTURE OF XIAO
}\\
\vskip 5pt
\end{center}
\vskip 2pt

\begin{center}
Miguel A. BARJA\footnotemark[1]\footnotetext[1]
{Partially supported by
CICYT PS93-0790 and HCM project n.ERBCHRXCT-940557.}.

{\scriptsize
Departament de Matem\`atica Aplicada I. Universitat Polit\`ecnica de
Catalunya. Barcelona. Spain}

\bigskip

Francesco ZUCCONI\footnotemark[2]\footnotetext[2]{Partially
supported by
HCM project n.ERBCHRXCT-940557.}.

{\scriptsize Dipartamento di Matematica e Informatica. Universit\`a degli Studi di
Udine. Udine. Italy}

\end{center}

\vskip 30pt
When $f:S\rightarrow B$ is a surjective morphism of a complex, smooth
surface $S$ onto a complex, smooth, genus $b$ curve $B$, such that the
fibre $F$ of $f$ has genus $g$, it is well known that
$f_{\ast}\omega_{S/B}={\cal E}$ is a locally free sheaf of
rank $g$ and degree $d= \X\O_{S}-(b-1)(g-1)$ and that $f$ is
not an holomorphic fibre bundle if and only if $d>0$.
In this case the {\it slope}, $\lambda (f)=\frac{K^{2}_{S}-8(b-1)(g-1)}
{d}$, is a natural invariant associated by Xiao to $f$
(cf. \cite{X1}). In \cite{X1}[Conjecture 2] he
conjectured that ${\cal{E}}$ has no locally free quotient
of degree zero (i.e., ${\cal E}$ is ample) if $\lambda(f)<4$. We give a partial affirmative answer to this conjecture:

\bigskip

{\bf Theorem 1.}\
{\it
Let $f:S\longrightarrow B$ be a relatively minimal
fibration with general fibre $F$. Let $b=g(B)$ and assume that $g=g(F) \geq 2$ and
that $f$ is not locally trivial.

If $\lambda (f) < 4$ then ${\cal E}=f_{\ast}\omega_{S/B}$ is ample provided one of the
following conditions hold
\begin{enumerate}
\item[(i)]  $F$ is non hyperelliptic.
\item[(ii)]  $b\leq 1$.
\item[(iii)]  $g(F)\leq 3$.
\end{enumerate}
}

\medskip

{\it Proof.} (i) If $q(S)>b$ the result follows from \cite{X1} Corollary
2.1. Now assume $q(S)=b$. By Fujita's decomposition
theorem (see \cite{F1}, \cite{F2} and also \cite{Ko} for a proof)

$${\cal E}={\cal A}\oplus{\cal F}_1\oplus ...\oplus {\cal F}_r$$
where $h^0(B,({\cal A}\oplus{\cal F}_1\oplus ...\oplus {\cal F}_r)
^{\ast})=0$,
${\cal A}$ is an ample sheaf and ${\cal F}_i$ are non trivial stable degree zero
sheaves. Then we only must prove that ${\cal F}_{i}=0$.
If $F$ is not hyperelliptic and rank $({\cal F}_i) \geq 2$
the claim is the content of \cite{X1} Proposition 3.1. If
rank$({\cal F}_i)=1$ we can use \cite{De} \S 4.2 or \cite{Ba}
Theorem 3.4 to conclude that ${\cal F}_i$ is torsion
in $\mbox{Pic}^0(B)$. Hence it induces an \'etale base change

$$
\xymatrix{
\widetilde{S}\ar[r]\ar[d]^{\widetilde{f}}&S\ar[d]^{f}\\
\widetilde{B}\ar[r]^{\sigma }&B}
$$

By flatness $\widetilde{f}_{\ast}\omega_{\widetilde{S}/\widetilde{B}}=\sigma
^{\ast}(f_{\ast }\omega_{S/B})$. Since $\sigma$ is \'etale $\lambda(f)=\lambda({\widetilde f})$ and $\sigma^{\ast}({\cal F}_i)={\cal O}_{\widetilde B}$ is a direct summand of ${\widetilde f}_{\ast}\omega_{{\widetilde S}/{\widetilde B}}$.
In particular by \cite{F1} $q({\widetilde S})>
{\widetilde b}=g({\widetilde B})$ hence $\lambda
({\widetilde f}) \geq 4$ by \cite{X1} Theorem 3.3: a contradiction.

(ii) If $b=0$  the claim is trivial. If $b=1$, any stable degree zero sheaf has rank one, then as in (i) we conclude.

(iii) If $g=2$ and ${\cal E}\neq{\cal A}$, then ${\cal E}={\cal A}\oplus\L$ where $\L$ torsion and we are done. The only non trivial case if $g=3$ is
${\cal E}={\cal A}\oplus\F$ where  $\A$ an ample line
bundle and  $\F$ a stable, degree zero, rank two vector bundle. Then $K^{2}_{S/B}\geq (2g-2)\mbox{deg }\A=4d$ and we are done by \cite{X1}[Theorem 2]\hfill $\Box$

\bigskip
\bigskip

Theorem 3.3 of \cite{X1} Xiao says that if $q(S)>b$ and
$\lambda(f)=4$ then ${\cal E}={\cal F}\oplus {\cal O}_B$, where
${\cal F}$ is a semistable sheaf. We have the following improvement:

\bigskip

{\bf Theorem 2.}\
{\it
Let $f:S \longrightarrow B$ be a relatively minimal non locally trivial
fibration. If $\lambda(f)=4$ then ${\cal E}=f_{\ast}\omega_{S/B}$ has
at most one degree zero, rank one quotient ${\cal L}$

Moreover, in this case ${\cal E}={\cal A}\oplus {\cal L}$ with
${\cal A}$ semistable and ${\cal L}$ torsion.
}

\medskip

{\it Proof.} As in the previous theorem the torsion subsheaf ${\cal L}$ becomes the trivial one after an \'etale base change; thus

$${\widetilde {f}}_{\ast}\omega_{{\widetilde S}/{\widetilde B}}=
{\widetilde {\cal A}}\oplus {\cal O}_{{\widetilde B}}, \quad
{\widetilde {\cal A}}=\sigma^{\ast}{\cal A}.$$

By \cite{X1}[Theorem 3.3], ${\widetilde {\cal A}}$ is semistable. Then ${\cal A}$ is also semistable by \cite{M1}[Proposition 3.2].
\hfill $\Box$

\end{document}